\documentstyle[12pt]{article}
\newtheorem{Lemma}{Lemma}
\newtheorem{Theorem}{Theorem}
\newtheorem{Prop}{Proposition}

\newtheorem{Rem}{Remark}
\newcommand{\pf}{\medskip\noindent{\sc Proof: }}
\newcommand{\qed}{$\Box$}
\newcommand{\rTo}{\,\longrightarrow\,}
\title{On Two Proofs for the Existence and Uniqueness of Integrals for Finite-Dimensional Hopf Algebras}
\author{Louis H. Kauffman\thanks{Research supported in part by NSF Grant
DMS 920-5227}  \phantom{aaa} and  \phantom{aaa}
David E. Radford\thanks{Research supported in part by NSF Grant
DMS 9802178}\\
\and \\ 
Department of Mathematics, Statistics \\
and Computer Science (m/c 249)    \\
851 South Morgan Street   \\
University of Illinois at Chicago\\
Chicago, Illinois 60607-7045}
\begin{document}
\maketitle
%
%\date{December 1, 1997}
Integrals play a basic role in the structure theory of finite-dimensional Hopf algebras $A$ and their duals $A^*$ over a field $k$. The existence and uniqueness of integrals for finite-dimensional Hopf algebras was first established by Hopf modules; for a long time the only means known for doing so. 

In
\cite{VDaele},
\cite{KuperRef}
and
\cite{NILL}
the existence and uniqueness of integrals for $A$ is established without using Hopf modules. We find the approach of 
\cite{KuperRef} 
very interesting in that it is based on a formalism which relates Hopf algebras and complete invariants of $3$-manifolds in a rather intriguing way.

This paper has two main purposes. The first is to explain the formalism of 
\cite{KuperRef}
to the extent that the proof of the existence and uniqueness of integrals for $A$ found in
\cite{KuperRef}
can be understood in more familiar algebraic terms. The reader is directed to
\cite{Fernando}
for a much fuller explanation of this formalism which includes a discussion of its subtleties and its connections with identities which hold in $A$. The description of integrals in 
\cite{KuperRef}
is given in terms of the trace function. The second purpose of this paper is to show that the theory of integrals for $A$ can be deduced from the ideas concerning the trace function on ${\rm End}(A)$ developed in
\cite[Section 5]{RTrace}.

Connections between Kuperberg's work and the trace function should be of interest to those who study $3$-manifold invariants which arise in some fashion from Hopf algebras. Treatments of integrals which do not involve Hopf modules should be of theoretical interest to those who study Hopf algebras. The ideas for this paper arose on examination of some of the remarkable connections between Hopf algebras and $3$-manifolds developed in 
\cite{KuperRef}.
The authors hope that the reader will see the advantage of using {\em all} of the formalisms discussed in this paper for the study of finite-dimensional Hopf algebras.

We shall assume that the reader has an elementary knowledge of Hopf algebras. A suggested reference is any one of 
\cite{Abe},
\cite{LR},
\cite{Mont}
and
\cite{SweedlerBook}.
Readers with a basic knowledge of Hopf algebras should find this paper virtually self contained. 

Throughout $k$ is a field, $A$ is a finite-dimensional Hopf algebra over $k$ and all vector spaces are over $k$.
\section{Preliminaries}\label{Preliminaries}
We write $M {\otimes}N$ for the tensor product $M {\otimes_k}N$ of vector spaces $M$ and $N$ and drop the subscript $k$ from ${\rm End}_k(M)$. The identity map of a vector space $M$ is denoted by $1_M$, which also denotes the unity of $M$ when $M$ is an algebra over $k$. The meaning of $1_M$ should always be clear from context.

For $a \in A$ we represent $\Delta (a) \in A {\otimes}A$ by $\Delta (a) = a_{(1)} {\otimes}a_{(2)}$. Observe that the dual algebra $A^*$ has a left $A$-module structure $(A^*, \rightharpoonup )$ and a right $A$-module structure $(A^*, \leftharpoonup )$ where
$$
(a{\rightharpoonup}p)(b) = p(ba) \quad \mbox{and} \quad (p{\leftharpoonup}a)(b) = p(ab)
$$
for all $a, b \in A$ and $p \in A^*$.  Likewise $A$ has a left $A^*$-module structure $(A, \rightharpoonup )$ and a right $A^*$-module structure $(A, \leftharpoonup )$ given by 
$$
p{\rightharpoonup}a = a_{(1)}p(a_{(2)}) \quad \mbox{and} \quad a{\leftharpoonup}p = p(a_{(1)})a_{(2)}
$$
for all $p \in A^*$ and $a \in A$. We let $\ell (a), r(a)$ be the endomorphisms of $A$ defined by 
$$
\ell (a)(b) = ab \quad \mbox{and} \quad r(a)(b) = ba
$$
for all $a, b \in A$ and we let $\ell (p), r(p)$ be the endomorphisms of $A$ defined by 
$$
\ell (p)(a) = p{\rightharpoonup}a \quad \mbox{and} \quad r(p)(a) = a{\leftharpoonup}p
$$
for all $p \in A^*$ and $a \in A$. Regard $A^* {\otimes}A$ as ${\rm End}(A)$ via the identification 
$$
(p {\otimes}a)(b) = p(b)a
$$
for all $p \in A^*$ and $a, b \in A$. Observe that ${\sf tr}(p {\otimes}a) = p(a)$ for $p{\otimes}a \in {\rm End}(A)$.

The vector space ${\rm End}(A)$ of linear endomorphisms of $A$ has a $k$-algebra structure, called the convolution algebra, whose product is defined by 
$$
(f{\star}g)(a) = f(a_{(1)})g(a_{(2)})
$$
for all $f, g \in {\rm End}(A)$ and $a \in A$. The unit of the convolution algebra is $\eta{\circ}\epsilon$, where $\eta : k \rTo A$ is the unit map which is defined by $\eta (1_k) = 1_A$.

Let $A^{op}$ denote the bialgebra over $k$ obtained from $A$ by replacing the product $m : A {\otimes}A \rTo A$ with $m^{op}$, where $a{\cdot}b = m^{op}(a {\otimes}b) = m(b {\otimes}a) = ba$ for all $a, b \in A$. Let $A^{cop}$ be the bialgebra over $k$ obtained from $A$ by replacing the coproduct $\Delta : A \rTo A {\otimes}A$ with $\Delta^{cop}$ defined by $\Delta^{cop}(a) = a_{(2)} {\otimes}a_{(1)}$ for all $a \in A$. Since $A$ is finite-dimensional the antipode $s$ of $A$ is bijective. Thus the bialgebras $A^{op}$ and $A^{cop}$ are Hopf algebras with antipode $s^{-1}$. Observe that $A^{op\, cop}$ is a Hopf algebra with antipode $s$ and $s : A \rTo A^{op\, cop}$ is a Hopf algebra map. The dual bialgebra $A^*$ is a Hopf algebra with antipode $S = s^*$.

An element $\Lambda \in A$ is a left (respectively right) integral for $A$ if $a\Lambda = \epsilon (a)\Lambda$ (respectively $\Lambda a = \Lambda\epsilon (a)$) for all $a \in A$. The set of left integrals for $A$ and the set of right integrals for $A$ are both ideals of $A$ which we denote by $\int_{\ell}$ and $\int_r$ respectively. We denote the set of left integrals for $A^*$ by $\int^{\ell}$ and we denote the set of right integrals for $A^*$ by $\int^r$. 

An example of a left integral described early on is Haar measure which is a functional. For this reason  left or right integrals for $A^*$ are sometimes referred to as left or right integrals, and left of right integrals for $A$ are referred to as left or right cointegrals. This is the convention of 
\cite{KuperRef}.

Observe that $\lambda \in \int^{\ell}$ if and only if $p\lambda = p(1)\lambda$ for all $p \in A^*$, or equivalently $a_{(1)}\lambda (a_{(2)}) = \lambda (a)1$ for all $a \in A$. Likewise $\lambda \in \int^r$ if and only if $\lambda p = (p(1))\lambda$ for all $p \in A^*$, or equivalently 
\begin{equation}\label{Eqlambda1}
\lambda (a_{(1)})a_{(2)} = \lambda (a)1 = \lambda (a) \eta (1_k)
\end{equation}
for all $a \in A$. Let $\lambda \in \int^r$. Since 
$$
\lambda (ab_{(1)})b_{(2)} = \lambda (a_{(1)}b_{(1)})s^{-1}(a_{(3)})a_{(2)}b_{(2)} = \lambda ((a_{(1)}b)_{(1)})s^{-1}(a_{(2)})(a_{(1)}b)_{(2)}
$$
for all $a, b \in A$, as a consequence of (\ref{Eqlambda1}) we have
\begin{equation}\label{EqlambdaabOne}
\lambda (ab_{(1)})b_{(2)} = \lambda (a_{(1)}b)s^{-1}(a_{(2)})
\end{equation}
for all $a, b \in A$.
\section{Integrals in Terms of the Trace Function}\label{SecIntTermsTrace}
Let $s$ be the antipode of $A$ and suppose that $\lambda \in \int^r$, $\Lambda \in \int_\ell$ are not zero. Then $\lambda (\Lambda ) \neq 0$ by 
\cite[Corollary 1]{RTrace}. 
We may assume that $\lambda (\Lambda ) = 1$. Part a) of 
\cite[Proposition 2]{RTrace} 
gives the equation
\begin{equation}\label{Eqs2Trace}
{\sf tr}(r(a){\circ} s^2 {\circ} r(p)) = \lambda (a) p(\Lambda )
\end{equation}
for all $a \in A$ and $p \in A^*$ which expresses a fundamental  connection between the trace function on ${\rm End}(A)$ and integrals. Motivated by ideas of 
\cite{KuperRef}
we express (\ref{Eqs2Trace}) in terms of the endomorphism $Q = \lambda {\otimes}\Lambda$ of $A$ as
\begin{equation}\label{Eqs2TraceQ}
p(Q(a)) = {\sf tr}(r(a){\circ} s^2 {\circ} r(p)) 
\end{equation}
for all $p \in A^*$ and $a \in A$ and we describe $Q$ in terms of structure constants. 

Let $\{ a_1, \ldots, a_n\}$ be a basis for $A$ and let $m^{\ell}_{i, j}, \Delta_{\ell}^{i, j}, s^i_j \in k$ for  $1 \leq i, j, \ell \leq n$ be the structure constants for the product $m$, the coproduct $\Delta$ and the antipode $s$ of $A$, defined by 
$$
m(a_\imath {\otimes}a_\jmath ) = a_\imath a_\jmath = \sum_{\ell = 1}^n m^{\ell}_{\imath, \jmath}a_{\ell}
$$
for all $1 \leq  \imath, \jmath \leq n$, 
$$
\Delta (a_{\ell}) = \sum_{\imath, \jmath =1}^n\Delta^{\imath, \jmath}_{\ell}a_\imath {\otimes}a_\jmath
$$
for all $1 \leq \ell \leq n$ and 
$$
s(a_\jmath ) = \sum_{\imath = 1}^ns^\imath_\jmath a_\imath
$$
for all $1 \leq \jmath \leq n$ respectively. Throughout this paper we shall use, whenever possible,  the Einstein summation convention which is to omit the summation symbol and to sum over the full range of the values of  those indices which occur as both upper and lower indices. Thus we write 
$$
a_\imath a_\jmath =  m^{\ell}_{\imath, \jmath}a_{\ell}, \quad \Delta (a_{\ell}) = \Delta^{\imath, \jmath}_{\ell}a_\imath {\otimes}a_\jmath \quad \mbox{and} \quad 
s(a_\jmath ) = s^\imath_\jmath a_\imath.
$$
Now let $\{\alpha_1, \ldots \alpha_n\}$ be the basis for $A^*$ dual to the basis $\{ a_1, \ldots, a_n\}$ for $A$. The reader can easily check that 
$$
(r(a_\jmath ){\circ} s^2 {\circ} r(\alpha_\imath ))(a_\ell ) = 
\Delta^{\imath, v}_\ell s^u_v s^w_u m^x_{w, \jmath}a_x
$$
for all $1 \leq \imath, \jmath, \ell  \leq n$. Therefore 
$$
{\sf tr}(r(a_\jmath ){\circ} s^2 {\circ} r(\alpha_\imath )) = 
m^\ell_{w, \jmath}\Delta^{\imath, v}_\ell s^u_v s^w_u
$$
for all $1 \leq \imath, \jmath \leq n$. Writing $Q(a_\jmath ) = Q^\imath_\jmath a_\imath$ observe that $Q^\imath_\jmath = \alpha_\imath (Q(a_\jmath ))$ for all $1 \leq \imath, \jmath \leq n$. Thus 
\begin{equation}\label{EqQijTr}
Q^\imath_\jmath = m^\ell_{w, \jmath}\Delta^{\imath, v}_\ell s^u_v s^w_u
\end{equation}
for all $1 \leq \imath, \jmath \leq n$.

Now let the decorated symbols
$$
\begin{array}{ccc}
\stackrel{i}{\searrow} & & \\
& m & \stackrel{\ell}\rightarrow \\
\stackrel{j}\nearrow & & 
\end{array}\;\;\; ,
\qquad
\begin{array}{ccc}
& & \stackrel{i}\nearrow \\
\stackrel{\ell}\rightarrow & \Delta & \\
& & \stackrel{j}\searrow
\end{array}
\qquad
\mbox{and} \;\;
\stackrel{j}\rightarrow s \stackrel{i}\rightarrow
$$
represent the structure constants $m^\ell_{\imath, \jmath}$, $\Delta^{\imath, \jmath}_\ell$ and $s^\jmath_\imath$ respectively. Then  (\ref{EqQijTr}) can be expressed symbolically as
$$
\raisebox{9ex}{$\stackrel{\jmath}{\rightarrow} Q \stackrel{\imath}{\rightarrow} \;\; = \;\;$}
\begin{picture}(140,100)(-70,-50)
\put(-47, -2){$m$}
\put(39, -2){$\Delta$}
\put(-2, 40){$s$}
\put(-2, -42){$s$}
\put(-61, -3){$\stackrel{\jmath}{\rightarrow}$}
\put(52, -3){$\stackrel{\imath}{\rightarrow}$}
\put(-5, 35){\vector(-1,-1){30}}
\put(8, -30){${}_{}^v$}
\put(35, -5){\vector(-1,-1){30}}
\put(-33, 10){${}_{}^w$}
\put(-35, 0){\vector(1,0){70}}
\put(24, -1){${}_{}^\ell$}
\put(0,-35){\vector(0,1){70}}
\put(2, 17){${}_{}^u$}
\end{picture} \;\; .
$$
%$$
%\stackrel{\jmath}\rightarrow Q \stackrel{\imath}\rightarrow \;\; = \;\;
%\begin{array}{ccccccc}
%                   & & & s & & & \\
%                   & & \stackrel{w}\swarrow & \stackrel{u}\uparrow & & & \\
%           \stackrel{\jmath}\rightarrow & m  &  &  & \stackrel{\ell}\rightarrow & \Delta & %\stackrel{\imath}\rightarrow \\
%                   & &  & & \stackrel{v}\swarrow & & \\
%                   & & & s & & & 
%\end{array} 
%$$
or, omitting the line labels,
$$
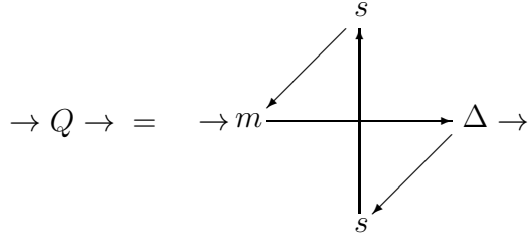
\raisebox{9ex}{$\rightarrow Q \rightarrow \;\; = \;\;$}
\begin{picture}(140,100)(-70,-50)
\put(-47, -2){$m$}
\put(39, -2){$\Delta$}
\put(-2, 40){$s$}
\put(-2, -42){$s$}
\put(-61, -3){$\rightarrow$}
\put(52, -3){$\rightarrow$}
\put(-5, 35){\vector(-1,-1){30}}
\put(35, -5){\vector(-1,-1){30}}
\put(-35, 0){\vector(1,0){70}}
\put(0,-35){\vector(0,1){70}}
\end{picture} \;\; .
$$
%\begin{array}{ccccccc}
%                   & & & s & & & \\
%                   & & \swarrow & \uparrow & & & \\
%          \rightarrow & m  & &  & \rightarrow & \Delta & \stackrel{\imath}\rightarrow \\
%                   & &  & & \swarrow & & \\
%                   & & & s & & & 
%\end{array} \;\;   .
%$$

Consider the endomorphism $P$ of $A$ defined by $P(a_\jmath ) = m^{\ell}_{w, \jmath}s^u_{\ell}\Delta^{v, x}_us^w_va_x$ for all $1 \leq \jmath \leq n$, or equivalently by
\begin{equation}\label{EqKupP}
P_\jmath^\imath = m^{\ell}_{w, \jmath}s^u_{\ell}\Delta^{v, \imath}_us^w_v
\end{equation}
for all $1 \leq \imath, \jmath \leq n$. Notice that $P$ can be described diagrammatically by 
$$
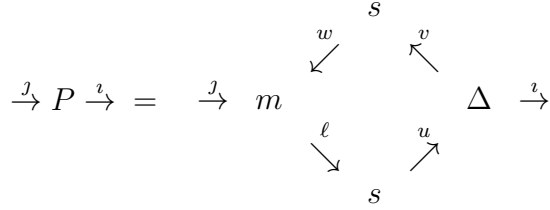
\stackrel{\jmath}\rightarrow P \stackrel{\imath}\rightarrow \;\; = \;\;
\begin{array}{ccccccc}
                   & & & s & & & \\
                   & & \stackrel{w}\swarrow & & \stackrel{v}\nwarrow & & \\
           \stackrel{\jmath}\rightarrow & m  & & & & \Delta & \stackrel{\imath}\rightarrow \\
                   & & \stackrel{\ell}\searrow & & \stackrel{u}\nearrow & & \\
                   & & & s & & & 
\end{array} 
$$
or, omitting the line labels, by
$$
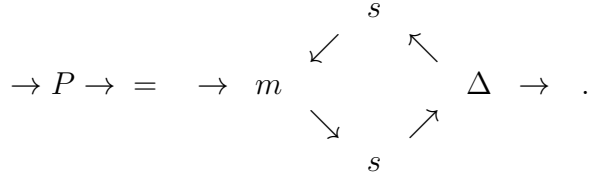
\rightarrow P \rightarrow \;\;  = \;\; \begin{array}{ccccccc}
                   & & & s & & & \\
                   & & \swarrow & & \nwarrow & & \\
           \rightarrow & m  & & & & \Delta & \rightarrow \\
                   & & \searrow & & \nearrow & & \\
                   & & & s & & & 
\end{array} \;\; .
$$
The diagrammatic descriptions of $Q$ and $P$, without the line labels, are examples of Kuperberg's formalism. By virtue of 
\cite[Lemma 3.3]{KuperRef} 
the endomorphism $P$ incorporates right integrals for $A^*$ and $A$, and furthermore ${\sf tr}P = 1$. It is shown more precisely that $P = \lambda {\otimes}\Lambda$, where $\lambda \in \int_r$, $\Lambda \in \int^r$ and $\lambda (\Lambda ) = 1$. In the next section we discuss the diagrammatic formalism of 
\cite{KuperRef} 
and Lemma \ref{LadderLem} which is the at the heart of the proof of existence and uniqueness of integrals given in 
\cite{KuperRef}.
\section{A Brief Discussion of Kuperberg's Diagrammatic Formalism for Hopf 
Objects}\label{SecKFormal}
In 
\cite{KuperRef}
Hopf objects are defined and their properties are developed. Many of the basic results for finite-dimensional Hopf algebras have analogs for Hopf objects. Here we briefly discuss the formalism for Hopf objects in the context of Hopf algebras and interpret it in terms of structure constants and symbolic computation with elements.

To begin, the structure maps of the finite-dimensional Hopf algebra $A$ are represented symbolically. The  product $m : A{\otimes A} \longrightarrow A$ and unit $\eta : A \longrightarrow k$, the coproduct $\Delta : A \longrightarrow A{\otimes} A$ and counit $\epsilon : A \longrightarrow k$, and the antipode $s : A \longrightarrow $A, are represented by
$$
\begin{array}{ccc}
\searrow & & \\ & m & \rightarrow \\ \nearrow & & 
\end{array}
\quad \mbox{and} \quad 
\begin{array}{cc} \eta & \rightarrow \end{array} \;\; ,
$$
$$
\begin{array}{ccc}
& & \nearrow  \\ & \rightarrow  \Delta &  \\ & & \searrow 
\end{array}
\quad \mbox{and} \quad 
\begin{array}{cc} \rightarrow \epsilon \end{array} \;\; ,
$$
and $\begin{array}{ccc} \rightarrow & s & \rightarrow \end{array}$ respectively. The incoming arrows can be thought of as inputs and the outgoing arrows can be thought of as outputs. Inputs for the product symbol are read {\em counterclockwise} and outputs for the coproduct symbol are read {\em clockwise}. An endomorphism $f$ of $A$ is represented $\begin{array}{ccc} \rightarrow & f & \rightarrow \end{array}$ and the identity of $A$ is represented by the shorthand notation $\longrightarrow$. If $f_1, \ldots, f_n$ are endomorphisms of $A$ then the endomorphism $f_1 {\otimes} \cdots {\otimes}f_n$ of $A^{\otimes n}$ is represented 
$$
\begin{array}{ccc}
\rightarrow & f_1 & \rightarrow \\
& \vdots & \\
\rightarrow & f_n & \rightarrow 
\end{array} \;\; .
$$

The diagrammatic representations of the structure maps of $A$ may be thought of in terms of structure constants, and combinations of these diagrams may be thought of as (sums of) products of structure constants. For example, 
$\begin{array}{ccc}
\searrow & & \\ & m & \rightarrow \\ \nearrow & & 
\end{array}$
may be thought of as $m^\ell_{\imath, \jmath}$ which is made explicit by writing 
$\begin{array}{ccc}
\stackrel{\imath}{\searrow} & & \\ & m & \stackrel{\ell}{\rightarrow} \\ \stackrel{\jmath}{\nearrow} & & 
\end{array}$. The associative law for multiplication $(a_\imath a_\jmath )a_k = a_\imath (a_\jmath a_k)$, which is equivalent to $m^\ell_{\imath, \jmath} m^u_{\ell, k}a_u = m^u_{\imath, \ell} m^\ell_{\jmath, k}a_u$, has the structure constant formulation $m^\ell_{\imath, \jmath} m^u_{\ell, k} = m^u_{\imath, \ell} m^\ell_{\jmath, k}$. The last equation is encoded in the diagram
$$
\begin{array}{ccccc}
\stackrel{\imath}{\searrow} & & & & \\
\stackrel{\jmath}{\rightarrow} & m & \stackrel{\ell}{\rightarrow} & m &  \stackrel{u}{\rightarrow}  \\
& &  \stackrel{k}{\nearrow} & &
\end{array}
\;\; = \;\;
\begin{array}{ccccc}
& & \stackrel{\imath}{\searrow}  & & \\
\stackrel{\jmath}{\rightarrow} & m & \stackrel{\ell}{\rightarrow} & m &  \stackrel{u}{\rightarrow}  \\
\stackrel{k}{\nearrow}& &  & &
\end{array} \;\; .
$$
Removing indices we have the diagrammatic formulation
$$
\begin{array}{ccccc}
\searrow & & & & \\
\rightarrow & m & \rightarrow & m & \rightarrow  \\
& &  \nearrow & &
\end{array}
\;\; = \;\;
\begin{array}{ccccc}
& & \searrow  & & \\
\rightarrow & m & \rightarrow & m &  \rightarrow  \\
\nearrow& &  & &
\end{array} 
$$
of the associative law. 

Two diagrams may be combined to produce a diagram by joining heads of (some) free output arrows of the first to tails of free input arrows of the second. The last equation illustrates this point. Diagrams with $m$ free input arrows and $n$ free output arrows may be regarded as morphisms from $A^{\otimes m}$ to $A^{\otimes n}$. Reading from top down, the $\imath^{th}$ free input (respectively output) arrow may be interpreted as the $\imath^{th}$ tensor factor on $A^{\otimes m}$ (respectively $A^{\otimes n}$).

A diagram may be modified to produce another diagram by joining heads of (some) free output arrows to tails of free input arrows. A very basic example. Let $f$ be an endomorphism of $A$. Since the structure constants of $f$ satisfy $f(a_\jmath ) = f^\imath_\jmath a_\imath$ it follows that ${\sf tr}f = f^\imath_\imath$, a scalar which may be regarded as the endomorphism of $A^{\otimes 0} = k$ defined by $1_k \mapsto {\sf tr}f$. The diagram 
$$
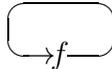
\begin{picture}(60,40)(-30,-5)
\put(-12.5, 15){\oval(15, 20)[l]}
\put(12.5, 15){\oval(15, 20)[r]}
\put(-12.5, 25){\line(1, 0){25}}
\put(2.5, 5){\line(1, 0){10}}
\put(-14, 2){$\rightarrow$}
\put(-4, 2){$f$}
\end{picture}
%\begin{picture}(40,40)(-20, -20)
%\put(0, 0){\circle{20}}
%\put(-5, -13){$f$}
%\end{picture}
$$
with no inputs and no outputs thought of in terms of structure constants represents ${\sf tr}f$.

The diagram 
$\begin{array}{ccc}
& & \nearrow  \\ & \rightarrow  \Delta &  \\ & & \searrow 
\end{array}$
can be thought of as $\Delta^\ell_{\imath, \jmath}$ which is made explicit by writing 
$\begin{array}{ccc}
& & \stackrel{\imath}{\nearrow}  \\ & \stackrel{\ell}{\rightarrow}  \Delta &  \\ & & \stackrel{\jmath}{\searrow}
\end{array}$.
It is not hard to see that the coassociative law can be expressed as 
$$
\begin{array}{ccccc}
& & \nearrow & & \\
\rightarrow & \Delta & \rightarrow & \Delta & \rightarrow \\
& & &  & \searrow
\end{array}
\;\;  =  \;\;
\begin{array}{ccccc}
& & & & \nearrow \\
\rightarrow & \Delta & \rightarrow & \Delta & \rightarrow \\
& & \searrow & &
\end{array}.
$$

The diagrams representing the structure maps for $A$ can also be thought of in terms of elements. Diagrams can be thought of as rules for symbolic computation. For example, 
$\begin{array}{ccc} \searrow & & \\ & m & \rightarrow \\ \nearrow & & \end{array}$ can be thought of as $m(a{\otimes}b) = c$ which we make explicit by writing 
$\begin{array}{ccc} \stackrel{a}{\searrow} & & \\ & m & \stackrel{ab}{\rightarrow} \\ \stackrel{b}{\nearrow} & & \end{array}$. With this interpretation the associative law
$$
\begin{array}{ccccc}
\stackrel{a}{\searrow} & & & & \\
\stackrel{b}{\rightarrow} & m & \stackrel{ab}{\rightarrow} & m &  \stackrel{(ab)c}{\rightarrow} \\
& & \stackrel{c}{\nearrow} 
\end{array}
\;\; = \;\;
\begin{array}{ccccc}
& & \stackrel{a}{\searrow} & & \\
\stackrel{b}{\rightarrow} & m & \stackrel{bc}{\rightarrow} & m &  \stackrel{a(bc)}{\rightarrow} \\
\stackrel{c}{\nearrow} & & & & 
\end{array}
$$
is more transparent.

The Einstein-like convention $\Delta (a) = a_{(1)}{\otimes}a_{(2)}$ for expressing the coproduct translates to 
$\begin{array}{ccc}
& & \stackrel{a_{(1)}}{\nearrow} \\
\stackrel{a}{\rightarrow} & \Delta & \\
& & \stackrel{a_{(2)}}{\searrow} 
\end{array}$.
The outputs in this case are symbolic elements. When thought of in terms of elements, combinations of diagrams {\em without directed closed loops} produce formal elements of tensor powers of $A$. As an example/exercise, the multiplicative property of the coproduct $\Delta (ab) = \Delta (a) \Delta (b)$ for all $a, b \in A$ is expressed by 
$$
\begin{array}{ccccc}
\searrow & & & & \nearrow \\
& m & \rightarrow & \Delta & \\
\nearrow & & & & \searrow 
\end{array} 
\;\; = \;\;
\begin{array}{ccccc}
\rightarrow & \Delta & \rightarrow & m & \rightarrow \\
& & \searrow \!\!\!\!\!\! \nearrow &  & \\
%& & \nearrow &  & \\
\rightarrow & \Delta & \rightarrow & m & \rightarrow 
\end{array} \quad , 
$$
where both diagrams are thought of in terms of elements or both diagrams are thought of in terms of structure constants. 

For diagrams without directed closed loops the basis-free way of thinking of them in terms of elements and the way of thinking of them in terms of structure constants are equivalent. For these diagrams it is useful to have both interpretations.

Using the canonical identification $A = A{\otimes} k$, we express the axiom $a1 = a$ for all $a \in A$, which is to say $m(1_A{\otimes}\eta ) = 1_A$, by
$$
\begin{array}{cccc}
& \searrow & & \\
\eta & \rightarrow & m & \rightarrow 
\end{array}
\;\; = \;\;
\rightarrow \quad  .
$$
The endomorphism $\eta {\circ} \epsilon$ of $A$ is expressed 
$\begin{array}{ccc} \rightarrow & \epsilon \;\; \eta & \rightarrow \end{array}$
and the endomorphism $\epsilon {\circ} \eta$ of $k$ is expressed as 
$\begin{array}{ccc} \eta & \rightarrow & \epsilon \end{array} = \eta \;\; \epsilon$. Using the notation conventions described above the reader is encouraged to formulate the other Hopf algebra axioms as diagrammatic equations. 

We conclude our discussion of Kuperberg's diagrammatic formalism at this point, perhaps having generated more questions than answers. The reader is directed to
\cite{KuperRef}
and 
\cite{Fernando}.

We end this section with a slight restatement of 
\cite[Lemma 3.1]{KuperRef},
which we refer to as the Ladder Lemma, for Hopf algebras. The Ladder Lemma is the basis for the proof of the existence and uniqueness of integrals described in 
\cite{KuperRef}.
\begin{Lemma}\label{LadderLem}
Each of the following four diagrams 
$$
\begin{array}{ccc} 
\rightarrow & m & \rightarrow \\ 
& \uparrow & \\
\rightarrow & \Delta & \rightarrow 
\end{array} \qquad
\begin{array}{ccc} 
\rightarrow & m & \rightarrow \\ 
& \uparrow & \\
\leftarrow & \Delta & \leftarrow 
\end{array}   \qquad
\begin{array}{ccc} 
\leftarrow & m & \leftarrow \\ 
& \uparrow & \\
\rightarrow & \Delta & \rightarrow 
\end{array} \qquad
\begin{array}{ccc} 
\leftarrow & m & \leftarrow \\ 
& \uparrow & \\
\leftarrow & \Delta & \leftarrow 
\end{array}
$$
is invertible, and in each case an inverse is obtained by replacing the vertical arrow $\uparrow$ with 
$\begin{array}{c} \uparrow \\ s \\ \uparrow \end{array}$ .
\end{Lemma}

\pf
We consider the first diagram, thinking of diagrams in terms of elements. We need to show that 
$$
\begin{array}{ccccc}
\rightarrow & m & \rightarrow & m & \rightarrow \\
& & & \uparrow & \\
& \uparrow & & s & \\
& & & \uparrow & \\
\rightarrow & \Delta & \rightarrow & \Delta & \rightarrow 
\end{array}
 \;\; = \;\;
\begin{array}{c} \rightarrow \\ \\  \\  \\ \rightarrow \end{array}
\;\; = \;\;
\begin{array}{ccccc}
\rightarrow & m & \rightarrow & m & \rightarrow \\
& \uparrow & & &  \\
& s & & \uparrow & \\
& \uparrow & & & \\
\rightarrow & \Delta & \rightarrow & \Delta & \rightarrow 
\end{array} \quad   .
$$
For inputs $a$ and $b$ the diagrammatic equations translate to 
$$
(ab_{(1)})s(b_{(2)(1)}){\otimes} b_{(2)(2)} = a{\otimes}b = (as(b_{(1)}))b_{(2)(1)} {\otimes} b_{(2)(2)},
$$
equations which follow directly from Hopf algebra axioms. The inverse statements for the other diagrams translate to 
$$
(ab_{(1)(2)})s(b_{(2)}){\otimes} b_{(1)(1)} = a{\otimes}b = (as(b_{(1)(2)}))b_{(2)} {\otimes} b_{(1)(1)},
$$
$$
b_{(1)}(s(b_{(2)(1)})a){\otimes} b_{(2)(2)} = a{\otimes}b = s(b_{(1)})(b_{(2)(1)}a) {\otimes} b_{(2)(2)}
$$
and 
$$
b_{(1)(2)}(s(b_{(2)})a){\otimes} b_{(1)(1)} = a{\otimes}b = s(b_{(1)(2)})(b_{(2)}a) {\otimes} b_{(1)(1)}
$$
respectively. These are easily verified as well.
\qed
\medskip

The diagrams and their inverses described in Lemma \ref{LadderLem} are called ladders in 
\cite{KuperRef}.
Note that the calculation 
\begin{eqnarray*}
(ab_{(1)})s(b_{(2)(1)}){\otimes} b_{(2)(2)} & = & a(b_{(1)}s(b_{(2)(1)})) {\otimes} b_{(2)(2)} \\
& = & a(b_{(1)(1)}s(b_{(1)(2)})) {\otimes} b_{(2)} \\
& = & a(\epsilon (b_{(1)})1) {\otimes} b_{(2)} \\
& = & a1 {\otimes}\epsilon (b_{(1)})b_{(2)} \\
& = & a {\otimes}\epsilon (b_{(1)})b_{(2)} \\
& = & a {\otimes}b
\end{eqnarray*}
is a blueprint for a diagrammatic verification of the first equation mentioned in the proof of Lemma \ref{LadderLem}.
\section{Integrals and Cointegrals}\label{SecIntCoInt}
Let $\lambda \in \int^r$, $\Lambda \in \int_r$ and consider the endomorphism ${\cal P} = \lambda {\otimes}\Lambda$ of $A$. Then
\begin{equation}\label{EqCalPIntCoInt}
{\cal P}(a_{(1)}) {\otimes}a_{(2)} = {\cal P}(a) {\otimes}1 \quad \mbox{and} \quad ({\cal P}(a))b = \epsilon (b){\cal P}(a)
\end{equation}
for all $a, b \in A$. The first equation follows by (\ref{Eqlambda1}) and the second follows by the definition of right integral for $A$. We follow the terminology of
\cite{KuperRef} 
and call any endomorphism of $A$ which satisfies (\ref{EqCalPIntCoInt}) an integral and cointegral for $A$.

It is easy to see that any linear combination of integrals and cointegrals for $A$ is again an integral and cointegral for $A$. Conversely:
\begin{Lemma}\label{LemCalP}
Let $A$ be a finite-dimensional Hopf algebra over the field $k$ and suppose that ${\cal P}$ is a non-zero integral and cointegral for $A$. Write ${\cal P} = \sum_{\imath = 1}^r \lambda_\imath {\otimes}\Lambda_\imath \in A^* {\otimes}A = {\rm End}(A)$ where $r$ is as small as possible. Then $\lambda_\imath \in \int^r$ and $\Lambda_\imath \in \int_r$ for all $1 \leq \imath \leq r$.
\end{Lemma}

\pf
Since ${\cal P} \neq 0$ and $r$ is as small as possible the sets $\{ \lambda_1, \ldots, \lambda_r\}$ and $\{ \Lambda_1, \ldots, \Lambda_r\}$ are linearly independent. Let $a \in A$. Then 
$$
\sum_{\imath = 1}^r \lambda_\imath (a_{(1)})\Lambda_\imath {\otimes}a_{(2)} = {\cal P}(a_{(1)}) {\otimes}a_{(2)} = {\cal P}(a) {\otimes}1 = \sum_{\imath = 1}^r \lambda_\imath (a)\Lambda_\imath {\otimes}1
$$
which means that 
$$
\sum_{\imath = 1}^r \Lambda_\imath {\otimes}\lambda_\imath (a_{(1)})a_{(2)} = \sum_{\imath =1}^r \Lambda_\imath {\otimes}\lambda_\imath (a)1.
$$
Since $\{ \Lambda_1, \ldots, \Lambda_r\}$ is linearly independent it follows that $\lambda_\imath (a_{(1)})a_{(2)} = \lambda_\imath (a)1$ for all $1 \leq \imath \leq r$. Thus $\lambda_\imath \in \int^r$ for all $1 \leq \imath \leq r$ by (\ref{Eqlambda1}).

Now let $b \in A$. Since $({\cal P}(a))b = \epsilon (b){\cal P}(a)$, or equivalently $\sum_{\imath = 1}^r \lambda_\imath (a)\Lambda_\imath b = \sum_{\imath = 1}^r \lambda_\imath (a)\Lambda_\imath \epsilon (b)$, for all $a \in A$, it follows that $\sum_{\imath = 1}^r \lambda_\imath {\otimes}\Lambda_\imath b = \sum_{\imath = 1}^r \lambda_\imath {\otimes}\epsilon (b)\Lambda_\imath$. Since $\{ \lambda_1, \ldots, \lambda_r \}$ is linearly independent we conclude that $\Lambda_\imath b = \epsilon (b)\Lambda_\imath$ for all $1 \leq \imath \leq r$. Thus $\Lambda_\imath \in \int_r$ for all $1 \leq \imath \leq r$.
\qed
\medskip

Notice that (\ref{Eqlambda1}) which defines $\lambda \in \int^r$ has the diagrammatic formulation 
$$
\begin{array}{cccc}
\rightarrow & \Delta & \rightarrow & \lambda \\
& & \searrow & \\
\end{array}
\;\; = \;\;
\rightarrow \lambda \;\; \eta \rightarrow.
$$
Let $\Lambda \in \int_r$. We will think of $\Lambda$ in terms of a linear function $\mu : k \longrightarrow A$ defined by $\mu (1_k) = \Lambda$. The fact that $\Lambda$ is a right integral for $A$ is expressed diagrammatically
$$
\begin{array}{cccc}
\mu & \rightarrow & m & \rightarrow \\
 & \nearrow & & \\
& & 
\end{array}
\;\; =  \;\;
\rightarrow \epsilon \;\; \mu \rightarrow.
$$
Notice that the equations of (\ref{EqCalPIntCoInt}) can be expressed by 
\begin{equation}\label{EqIntCoIntDiagDef}
\begin{array}{ccc}
\leftarrow & P & \\
& \uparrow &  \\
\rightarrow & \Delta & \rightarrow
\end{array}
\; = \;
\begin{array}{cccc}
\nwarrow & & & \\
& P & \eta & \rightarrow \\
\nearrow & & &
\end{array}
\;\;  \mbox{and} \;\;
\begin{array}{ccc}
\rightarrow & P & \\
& \downarrow  & \\
\rightarrow & m & \rightarrow
\end{array}
\; = \;
\begin{array}{cccc}
& & \searrow & \\
\rightarrow & \epsilon & & P \\
& & \swarrow  & 
\end{array} \;\;   .
\end{equation}
We are now ready to state and prove the existence and uniqueness for integrals and cointegrals along the lines of
\cite{KuperRef}. 
Part a) of the next theorem is 
\cite[Lemma 2.3]{KuperRef}
and part b) follows by 
\cite[Corollary 3.5]{KuperRef}. 
\begin{Theorem}\label{ThmKuP}
Let $A$ be a finite-dimensional Hopf algebra and $P$ be the endomorphism of defined by (\ref{EqKupP}). Then:
\begin{enumerate}
\item[{\rm a)}]
$P$ is an integral and cointegral for $A$ which satisfies ${\sf tr}P = 1$.
\item[{\rm b)}]
Integrals and cointegrals for $A$ are unique up to scalar multiple.
\end{enumerate}
\end{Theorem}

\pf
We follow the proof given in 
\cite{KuperRef}
and spell out a few more details in the process. The following arguments can be viewed as diagrammatic blueprints for structure constant proofs, or they can be viewed as a valid symbolic proofs in their own right. 

To show part a), we begin by establishing the first equation of (\ref{EqIntCoIntDiagDef}). It suffices to show that 
$$
\begin{array}{ccccc}
\leftarrow & \Delta & \leftarrow & P & \\
& \downarrow & & \uparrow & \\
\rightarrow & m & \rightarrow & \Delta & \rightarrow 
\end{array}
\;\; = \;\;
\begin{array}{cccc}
\nwarrow &  & &  \\
&  s & \eta & \rightarrow  \\
\nearrow & & & 
\end{array} \;\;    .
$$
For multiplying both sides of this equation on the left by 
$\begin{array}{ccc}
\leftarrow & \Delta & \leftarrow \\
& \downarrow & \\
& s & \\
& \downarrow & \\
\rightarrow & m &\rightarrow
\end{array} \,\, ,
$
which is an inverse of 
$\begin{array}{ccc}
\leftarrow & \Delta & \leftarrow \\
& \downarrow & \\
\rightarrow & m &\rightarrow
\end{array} 
$
by Lemma \ref{LadderLem}, we obtain the first equation of (\ref{EqIntCoIntDiagDef}). We calculate
\begin{eqnarray*}
\begin{array}{ccccc}
\leftarrow & \Delta & \leftarrow & P & \\
& \downarrow & & \uparrow & \\
\rightarrow & m & \rightarrow & \Delta & \rightarrow 
\end{array}
& = & 
\begin{array}{ccccccc}
& & & \Delta & \leftarrow & s & \\
& & \swarrow & \downarrow & & \uparrow & \\
\leftarrow & \Delta & & s & \rightarrow & m & \\
& \downarrow & & & & \uparrow & \\
\rightarrow & m & & \rightarrow & & \Delta & \rightarrow 
\end{array} \\ \\  \\ 
& = & 
\begin{array}{ccccccc}
& & & \Delta & \leftarrow & s & \\
& & \swarrow & \downarrow & & \uparrow & \\
\leftarrow & \Delta & & s & \rightarrow & m & \\
& \downarrow & & & & \uparrow & \\
& \Delta & & \rightarrow & & m & \rightarrow  \\
& & & \searrow \!\!\!\!\!\! \nearrow & &  & \\
\rightarrow & \Delta & & \rightarrow & & m & \rightarrow  
\end{array} \\ \\  \\ 
& = & 
\begin{array}{ccccccc}
& & \leftarrow & \Delta & \leftarrow & s & \\
& & \swarrow &  & & \uparrow & \\
& \Delta & \rightarrow & s & \rightarrow & m & \\
& \downarrow & & & & \uparrow & \\
& \Delta & & \rightarrow & & m & \rightarrow  \\
& & & \searrow \!\!\!\!\!\! \nearrow & &  & \\
\rightarrow & \Delta & & \rightarrow & & m & \rightarrow  
\end{array} \\ \\  \\ 
& = & \begin{array}{ccccc}
\leftarrow & \Delta & \leftarrow & s & \\
& & \searrow \!\!\!\!\!\! \nearrow &  & \\
\rightarrow & \Delta & \rightarrow & m & \rightarrow
\end{array} \\ \\  \\ 
& = & 
\begin{array}{ccccc}
\nwarrow & & & & \\
& s &  \eta & \rightarrow & \\
\nearrow & & & &  \quad .
\end{array}
\end{eqnarray*}
The first equation of our calculation follows since the coproduct is multiplicative, the second since the coproduct is coassociative and the third by Lemma \ref{LadderLem}. The diagrams in the last equation have no closed directed loops and the equation is $s(a_{(1)})_{(2)}{\otimes}s(a_{(1)})_{(1)}a_{(2)} = s(a){\otimes}1$ when thought of in terms of elements. The reader is encouraged to justify the last equation in our calculation using diagrammatic equivalents of Hopf algebra axioms.

Our calculation establishes the first equation of (\ref{EqIntCoIntDiagDef}). The reader is left with the task of establishing the second.

Since 
$\begin{array}{ccccccc}
& & & s & & & \\
& & \nearrow & & \searrow & & \\
\rightarrow & \Delta & & \rightarrow & & m & \rightarrow
\end{array} = 
\begin{array}{ccc}
\rightarrow & \epsilon \;\;  \eta & \rightarrow
\end{array}$,
the diagrammatic equivalent of the Hopf algebra axiom $s(a_{(1)})a_{(2)} = \epsilon (a)1$, we calculate 
\medskip

\begin{eqnarray*}
{\sf tr} P & = &  \raisebox{-6ex}{
\begin{picture}(150,80)(-75,-30)
\put(-47, -2){$m$}
\put(38, -2){$\Delta$}
\put(-2, 40){$s$}
\put(-2, -42){$s$}
\put(-60, 0){$\vector(1,0){10}$}
\put(-60, 0){$\vector(1,0){10}$}
\put(-60, 10){\oval(20,20)[bl]}
\put(-70, 10){\line(0,1){30}}
\put(-60, 40){\oval(20,20)[tl]}
\put(-60, 50){\line(1,0){120}}
\put(50, 0){$\line(1,0){10}$}
\put(60, 10){\oval(20,20)[br]}
\put(70, 10){\line(0,1){30}}
\put(60, 40){\oval(20,20)[tr]}
\put(-5, 35){\vector(-1,-1){30}}
\put(35, 5){\vector(-1,1){30}}
\put(-35, -5){\vector(1,-1){30}}
\put(5, -35){\vector(1,1){30}}
%\put(-35, 0){\vector(1,0){70}}
%\put(0,-35){\vector(0,1){70}}
\end{picture} } \\  
& = & 
\begin{picture}(100,100)(-50,-50)
\put(-47, -2){$m$}
\put(39, -2){$\Delta$}
\put(-3, -88){$s$}
\put(-3, -42){$s$}
\put(35, -5){\vector(-1,-1){30}}
%\put(-35, -5){\vector(1,-1){30}}
\put(35, 0){\vector(-1,0){62}}
\put(-5, -35){\vector(-1,1){30}}
\put(5, -80){\vector(1,2){36}}
\put(-41, -8){\vector(1,-2){36}}
\end{picture} \;\;  \\ \\ \\ \\ 
& = & 
\begin{array}{ccccc}
\eta & & & & \epsilon \\
& \searrow & & \swarrow &\\
& & s & &
\end{array} \\ \\ \\ 
& = & 
\eta \rightarrow s \rightarrow \epsilon \\ \\ 
& = & 
\eta \rightarrow \epsilon
\end{eqnarray*}
which establishes ${\sf tr} P = 1$. We have shown part a).

To prove part b) we first note that 
$$
\begin{array}{cccc}
\Lambda & \leftarrow & m & \leftarrow \\
& & \uparrow & \\
\mu & \rightarrow & \Delta & \rightarrow 
\end{array} =
\begin{array}{cccccc}
& & & & & \swarrow \\
\mu & \rightarrow & \lambda  & & s & \\
& & & & & \searrow 
\end{array} \quad ,
$$
which is the diagrammatic expression for the Hopf algebra identity derived from 
\begin{eqnarray*}
\lambda (\Lambda_{(1)}a)\Lambda_{(2)} & = &
\lambda (\Lambda_{(1)}a_{(1)})\Lambda_{(2)}a_{(2)}s(a_{(3)}) \\
& = & 
\lambda ((\Lambda a_{(1)})_{(1)})(\Lambda a_{(1)})_{(2)}s(a_{(2)}) \\
& = & 
\lambda (\epsilon (a_{(1)})\Lambda_{(1)})\Lambda_{(2)}s(a_{(2)}) \\
& =& 
\lambda (\Lambda_{(1)})\Lambda_{(2)}s(a) \\
& =& 
\lambda (\Lambda )s(a),
\end{eqnarray*}
a fundamental calculation which shows how the antipode is expressed in terms of cointegrals and integrals when $\lambda (\Lambda ) = 1$. Note that this calculation involves a ladder and its inverse; thus a diagrammatic proof of the Hopf algebra identity naturally involves Lemma \ref{LadderLem}. Multiplying both sides of the expression on the right by 
$\begin{array}{ccc}
\leftarrow & m & \leftarrow \\
& \uparrow & \\
& s & \\
& \uparrow & \\
\rightarrow & \Delta & \rightarrow
\end{array}$,
the inverse of 
$\begin{array}{ccc}
\leftarrow & m & \leftarrow \\
& \uparrow & \\
\rightarrow & \Delta & \rightarrow
\end{array}$,
gives the equation
\medskip

$$
\begin{array}{cc}
\lambda & \leftarrow \\
\mu & \rightarrow
\end{array} 
=
\begin{array}{ccccc}
\mu & \rightarrow & \lambda & & P
\end{array} \quad .
$$
We have established part b) and the theorem is proved.
\qed
\begin{Rem}
In proving the uniqueness of $P$ we effectively give a clear motivation for the diagrammatic formula which describes $P$.
\end{Rem}

For given a right integral $\lambda$ and a right cointegral $\mu$ we calculate that 
$$
\begin{array}{ccccccc}
\mu & \rightarrow & \Delta           & \rightarrow & m             & \rightarrow  & \lambda \\
      &                   & \downarrow &                    & \uparrow &                    &              
\end{array}
\;\; = \;\;
\begin{array}{ccc}
\mu & \rightarrow & \lambda \\
\rightarrow & s & \rightarrow
\end{array}.
$$
Applying the appropriate inverse ladder we obtain
\begin{eqnarray*}
\begin{array}{cc}
\mu & \rightarrow \\ \rightarrow & \lambda 
\end{array}
& = &
\begin{array}{ccccccccc}
\mu & \rightarrow & \Delta           & &  \rightarrow &  & m             & \rightarrow  & \lambda \\
      &                   & \downarrow &                   &  & & \uparrow                    &      &        \\
      & \leftarrow   & \Delta           & \rightarrow & s & \rightarrow & m & \leftarrow &
\end{array} \\ \\ \\
& = & 
\begin{array}{ccccccc}
& & \mu & \rightarrow & \lambda & & \\
& & & s & & &  \\
& & \swarrow & & \nwarrow & & \\
\leftarrow & \Delta & \rightarrow & s & \rightarrow & m & \leftarrow 
\end{array} \\ \\  \\ 
& = & 
\begin{array}{ccc}
\mu & \rightarrow & \lambda \\
\leftarrow & P & \leftarrow 
\end{array}.
\end{eqnarray*}
Thus the structure constant formula for $P$ appears naturally when one considers this general property of integrals and cointegrals.
\begin{Rem}
The endomorphism $P$ of $A$ can be described in terms of an endomorphism ${\cal Q}$ of $A{\otimes}A$ and the trace function.
\end{Rem}

This is done as follows. Cut the arrow 
$\begin{array}{ccc} & & s \\ & \swarrow & \\ m & & \end{array}$ in the unlabeled diagram found after (\ref{EqKupP}) which defines $P$ to obtain an endomorphism ${\cal Q}$ of $A{\otimes}A$ determined by 
\begin{eqnarray*}
\begin{array}{ccc}
\searrow & & \nearrow \\ 
& {\cal Q} & \\
\nearrow  & & \searrow
\end{array}
& = &
\raisebox{-9ex}{
\begin{picture}(140,100)(-70,-50)
\put(-47, -2){$m$}
\put(39, -2){$\Delta$}
\put(-2, 40){$s$}
\put(-2, -42){$s$}
\put(-61, -3){$\rightarrow$}
\put(52, -3){$\rightarrow$}
\put(-5, 35){\vector(-1,-1){12}}
\put(-22, 18){\vector(-1,-1){12}}
\put(35, 6){\vector(-1,1){30}}
\put(-35, -3){\vector(1,-1){30}}
\put(5,-32){\vector(1,1){30}}
\end{picture}} \\  
%\raisebox{9.2ex}{= $\;\;$}
& = & 
\raisebox{1.2ex}{
$\begin{array}{ccccccccc}
& & & & & &  & s &  \rightarrow \\ 
\searrow & & & & & &  \nearrow & &  \\ 
& m & \rightarrow & s & \rightarrow & \Delta & & & \\
\nearrow & & & & & & \searrow & &  \\ 
\end{array}.$}
\end{eqnarray*}
Thus ${\cal Q}(a{\otimes}b) = s(s(ab)_{(1)}) {\otimes} s(ab)_{(2)}$ for all $a, b \in A$. It is an easy exercise to show that 
$$
P = ({\sf tr}{\otimes} 1_{{\rm End}(A)})({\cal Q}),
$$
where we view ${\sf tr}{\otimes} 1_{{\rm End}(A)}$ as the composite 
$$
{\rm End}(A{\otimes}A) = {\rm End}(A){\otimes}{\rm End}(A) \stackrel{{\sf tr} {\otimes} 1_{{\rm End}(A)}}{\rTo} {\rm End}(A).
$$
This calculation suggests how closed directed loops in diagrams can be realized as composites involving the trace function and maps of tensor powers of $A$. In general, this method of cutting loops and introducing traces provides a powerful method for translation between the tensor diagrammatic and pure algebraic approaches to Hopf algebras. By this method one can translate a diagrammatic proof into a purely algebraic proof by cutting closed directed loops in the diagram and introducing traces to the algebraic formalism.

The diagrammatic approach to Hopf algebras is formulated in category theoretic terms in
\cite{KuperRef}. In this framework each diagram is a morphism in the category with input arrows emanating from a domain object and output arrows terminating in a range object. The category is an associative tensor category generated by objects $V$ and $k$ which can be thought of as a representation space for $A$ (and therefore $A$ itself) and the field of definition for $A$. The axioms imply identifications $k {\otimes}V = V = V{\otimes}k$ and $k{\otimes}k = k$. 

A diagram with no free input arrows is a morphism from $k$. A diagram with no free output arrows is a morphism to $k$. Thus a diagram with no free input or output arrows is a morphism from $k$ to $k$. A diagram with $n$ free inputs arrow and $m$ free output arrows is a morphism from $V^{\otimes n}$ to $V^{\otimes m}$.

If we work in the full category of diagrams with closed directed loops then the functors to algebraic categories are primarily restricted to categories associated with finite-dimensional Hopf algebras since the functors must be defined via the structure constants (tensors) of the algebra. There are many interesting questions which can be asked about the faithfulness of these functors.
\section{Another Existence and Uniqueness Proof for Integrals}\label{SecAnotIntCoInt}
As we have mentioned, the theory of integrals for finite-dimensional Hopf algebras over the field $k$ was originally treated in the context of Hopf modules
\cite{LS}, \cite{SweedlerBook}. 
The theory of Hopf modules accounts for a crucial result, which can be rephrased as 
\cite[Theorem 1]{RTrace},
from which the general theory of integrals is derived. The theorem we are citing is that there exists a linear isomorphism $f : A \rTo A^*$ which satisfies
\begin{equation}\label{EqAIsoAStar}
f(ab) = f(b){\leftharpoonup}s(a) \quad \mbox{and} \quad f(a{\leftharpoonup}p) = f(a)p
\end{equation}
for all $a, b \in A$ and $p \in A^*$. As a consequence $s$ is bijective.

We now set the stage for an algebraic proof of 
\cite[Theorem 1]{RTrace}
which uses the trace function and does not involve Hopf modules. We first show that $s$ is bijective using a convolution algebra argument. See 
\cite[Section 3]{RO}
and 
\cite[Theorem 1.6.2]{LR}.
Our argument paraphrases the proof given for 
\cite[Theorem 1.6.2]{LR}.
\begin{Lemma}\label{SIsBijective}
Let $A$ be a finite-dimensional Hopf algebra with antipode $s$ over the field $k$. Then $s$ is bijective.
\end{Lemma}

\pf
We prove the lemma by induction on ${\rm Dim}\,A$. The case ${\rm Dim}\,A = 1$ is trivial.

Suppose that all Hopf algebras whose dimension is less than ${\rm Dim}\,A$ have bijective antipodes. Consider the Hopf algebra $B = s(A)$ with antipode $s|_B$. If $B = A$ then $s = s|_B$ and is onto; hence $s_B$ is bijective. If $B \neq A$ then $s|_B : B \rTo B$ is bijective by our induction hypothesis. In either case $A = {\rm Ker}s \oplus B$ as vector spaces.

Let $\pi : A \rTo B$ be the linear projection onto $B$ with kernel ${\rm Ker}s$. Since $\pi, s$ and $\epsilon$ vanish on ${\rm Ker}s$ and $\Delta ({\rm Ker}s) \subseteq {\rm Ker}s {\otimes}A + A \otimes{\rm Ker}s$ we conclude that $\pi{\star}s \equiv \eta{\circ}\epsilon$ on ${\rm Ker}s$. Since $\pi \equiv 1_A$ on $B$ and $\Delta (B) \subseteq B {\otimes}B$ we conclude that $\pi{\star}s \equiv  \eta{\circ}\epsilon$ on $B$. Therefore  $\pi{\star}s = \eta{\circ}\epsilon$ which means that $\pi$ is a left inverse for $s$ in the convolution algebra ${\rm End}(A)$. Hence $\pi = 1_A$, or equivalently ${\rm Ker}s = (0)$. We have shown that $s$ is one-one. This implies that $s$ is bijective as required.
\qed
\medskip

Part a) of the next result has its origin in
\cite{LS}.
\begin{Prop}\label{LambdaLambdaIs1}
Let $A$ be a finite-dimensional Hopf algebra over the field $k$. Suppose that $\lambda \in \int^r$ and $\Lambda \in \int_r$ satisfy $\lambda (\Lambda ) = 1$. Then:
\begin{enumerate}
\item[{\rm a)}]
$s(a) = \Lambda {\leftharpoonup}(a{\rightharpoonup}\lambda )$ for all $a \in A$.
\item[{\rm b)}]
The endomorphism $f$ of $A$ defined by $f(a) = \lambda{\leftharpoonup}s(a)$ for all $a \in A$ satisfies (\ref{EqAIsoAStar}).
\end{enumerate}
\end{Prop}

\pf
Since $\lambda (\Lambda ) =1$, part a) follows by the calculation made in the proof of part b) of Theorem \ref{ThmKuP}. As for part b), observe that the first equation of (\ref{EqAIsoAStar}) follows since $s : A \longrightarrow A^{op}$ is an algebra map. We use (\ref{EqlambdaabOne}) and the fact that $s : A \rTo A^{op\, cop}$ is a bialgebra map to establish the second equation of (\ref{EqAIsoAStar}) in the calculation
\begin{eqnarray*}
(f(a)p)(b)  & = & f(a)(b_{(1)})p(b_{(2)})  \\
& = & \lambda (s(a)b_{(1)})p(b_{(2)})  \\
& = & \lambda (s(a)_{(1)}b)p(s^{-1}(s(a)_{(2)}))  \\
& = & \lambda (s(a_{(2)})b)p(s^{-1}(s(a_{(1)})))  \\
& = & \lambda (s(a{\leftharpoonup}p)b)  \\
& = & f(a{\leftharpoonup}p)(b)
\end{eqnarray*}
for all $a, b \in A$ and $p \in A^*$.

It remains to show that $f$ is bijective, or equivalently that $f$ is one-one. Define an endomorphism $g$ of $A$ by $g(a) = a{\rightharpoonup}\lambda$ for all $a \in A$. Since $s$ is bijective, by part a) we conclude that $g$ is one-one and consequently is bijective. Since $f(a)(b) = \lambda (s(a)b) = g(b)(s(a))$ for all $a, b \in A$ it now follows that $f$ is one-one and is therefore bijective.
\qed
\medskip

Let ${\cal E} = {\cal E}_A$ be the endomorphism of ${\rm End}(A)$ defined by 
$$
p({\cal E}(f)(a)) = {\sf tr}(\ell (a){\circ}f{\circ} r(p))
$$
for all $p \in A^*$, $f \in {\rm End}(A)$ and $a \in A$. The endomorphism $E = E_H$ defined and studied in 
\cite[Section 5]{RTrace}
is ${\cal E}_{A^{op}}$.

In 
\cite[Section 5]{RTrace}
a detailed study of the endomorphism $E$ is made based on properties of integrals. Here we take the opposite tact. We use different arguments to establish basic properties of ${\cal E}$ and use these properties to deduce the existence of right integrals $\lambda \in \int^r$ and $\Lambda \in \int_r$ which satisfy $\lambda (\Lambda ) = 1$.

Let $a \in A$ and $p \in A^*$. In order to derive basic properties of ${\cal E}$ we consider the relationship between the endomorphisms $r(a)$ and $\ell (p)$ of $A$ and define right $A$-module and $A^*$-module actions on ${\rm End}(A)$. Note that
\begin{equation}\label{EqRaEllP}
r(a){\circ}\ell (p) = 1_A\star (p {\otimes}a) =  \ell (s(a_{(2)}){\rightharpoonup}p){\circ}r(a_{(1)}).
\end{equation}
Now let $f \in {\rm End}(A)$. Using the fact that $r(a){\circ}r(b) = r(ba)$ for all $b  \in A$ and the fact that $s : A \rTo A^{op}$ is an algebra map it is easy to see that ${\rm End}(A)$ is a right $A$-module according to 
$$
f{\bullet}a = r(a_{(1)}){\circ}f{\circ}r(s(a_{(2)})).
$$
Likewise ${\rm End}(A)$ is a right $A^*$-module according to 
$$
f{\bullet}p = \ell (S(p_{(2)})){\circ}f{\circ}\ell (p_{(1)}).
$$
For the proof of part c) of the next proposition we will use the relation
\begin{equation}\label{EqRR}
r(p{\circ}r(a)) = r(s(a_{(2)})){\circ} r(p){\circ}r(a_{(1)}).
\end{equation}
\begin{Prop}\label{EProperties}
Let $A$ be a finite-dimensional Hopf algebra over the field $k$. Then:
\begin{enumerate}
\item[{\rm a)}]
${\cal E}(p {\otimes}a) = \ell (p){\circ}r(a)$ for all $p \in A^*$ and $a \in A$.
\item[{\rm b)}]
${\cal E}$ is a linear automorphism of ${\rm End}(A)$.
\item[{\rm c)}]
$r(a){\circ}{\cal E}(f) = {\cal E}(f{\bullet}a)$ for all $a \in A$ and $f \in {\rm End}(A)$.
\item[{\rm d)}]
${\cal E}(f){\circ}\ell (p) = {\cal E}(f{\bullet}p)$ for all $f \in {\rm End}(A)$ and $p \in A^*$.
\end{enumerate}
\end{Prop}

\pf
Let $p, q \in A^*$ and $a, b \in A$. We first show part a). Using the fact that 
$$
\ell (b){\circ}(p {\otimes}a){\circ}r(q) = qp {\otimes}ba
$$
we compute 
%\begin{eqnarray*}
$$
q({\cal E}(p {\otimes}a)(b))  =  {\sf tr}(\ell (b){\circ} (p {\otimes}a){\circ}r(q)) 
 = {\sf tr}(qp {\otimes}ba)  = qp(ba).
$$
%\end{eqnarray*}
On the other hand 
%\begin{eqnarray*}
$$
q((\ell (p){\circ}r(a))(b))  =  q(\ell (p)(ba)) 
 =  q(p{\rightharpoonup}(ba)) 
 =  qp(ba).
$$
%\end{eqnarray*}
Thus $q({\cal E}(p {\otimes}a)(b)) = q((\ell (p){\circ} r(a))(b))$ which establishes part a).
To show part b) we note that the identity map $1_A$ has a left inverse in the convolution algebra ${\rm End}(A)$, namely the antipode $s$ of $A$. Therefore by (\ref{EqRaEllP}) the $r(a){\circ}\ell (p)$'s and also the $\ell (p) {\circ} r(a)$'s span ${\rm End}(A)$; in particular the $\ell (p){\circ}r(a)$'s span ${\rm End}(A)$. Thus by part a) the endomorphism ${\cal E}$ of ${\rm End} (A)$ is onto.  Since $A$ is finite-dimensional ${\cal E}$ must be a linear automorphism of ${\rm End}(A)$. Part b) is established.

Let $f \in {\rm End}(A)$ and $a \in A$. We use (\ref{EqRR}) and the fact that $\ell (b)$ and $r(c)$ commute for all $b, c \in A$ to compute
\begin{eqnarray*}
q((r(a){\circ}{\cal E}(f))(b)) & = & (q{\circ}r(a))({\cal E}(f)(b)) \\
& = & {\sf tr}(\ell (b){\circ}f{\circ} r(q{\circ}r(a))) \\
& = & {\sf tr}(\ell (b){\circ}f{\circ}r(s(a_{(2)})){\circ} r(q){\circ}r(a_{(1)})) \\
& = & {\sf tr}(r(a_{(1)}){\circ}\ell (b){\circ}f{\circ}r(s(a_{(2)})){\circ} r(q)) \\
& = & {\sf tr}(\ell (b){\circ}r(a_{(1)}){\circ}f{\circ}r(s(a_{(2)})){\circ} r(q)) \\
& = & q({\cal E}(f{\bullet}a)(b))
\end{eqnarray*}
for all $q \in A^*$ and $b \in A$ which establishes part c). Part d) follows in a similar manner. Better yet, observe that part d) follows from part c) applied to $A^*$. 
\qed
\medskip

We next observe that
\begin{equation}\label{EqSMinus2One}
s^{-2}{\bullet}a = \epsilon (a)s^{-2} \quad \mbox{and} \quad 
s^{-2}{\bullet}p = p(1)s^{-2} 
\end{equation}
for all $a \in A$ and $p \in A^*$. To establish the first equation of (\ref{EqSMinus2One}) let $a \in A$. Since $s^{-2}$ is an algebra endomorphism of $A$ and $s^{-1}$ is an antipode for $A^{op}$ it follows that 
\begin{eqnarray*}
(s^{-2}{\bullet}a)(b) & = & (r(a_{(1)}){\circ} s^{-2} {\circ} r(s(a_{(2)})) )(b) \\
& = & s^{-2}(bs(a_{(2)}))a_{(1)} \\
& = & s^{-2}(b)s^{-1}(a_{(2)})a_{(1)} \\
& = & s^{-2}(b)(\epsilon (a)1) \\
& = & (\epsilon (a)s^{-2})(b)
\end{eqnarray*}
for all $b \in A$ which implies $s^{-2}{\bullet}a = \epsilon (a)s^{-2}$.  The reader is left with the exercise of establishing the second equation of (\ref{EqSMinus2One}).

Let ${\cal P} = {\cal E}(s^{-2})$. Using (\ref{EqSMinus2One}) and parts c), d) of Proposition \ref{EProperties} we see that $r(a){\circ} {\cal P} = \epsilon (a){\cal P}$ for all $a \in A$ and ${\cal P}{\circ}\ell (p) = p(1){\cal P}$ for all $p \in A^*$. Thus ${\cal P}$ is an integral and cointegral for $A$. Since ${\cal E}_{A^{op\, cop}}({\cal P}) \neq 0$ and
$$
p({\cal E}_{A^{op\, cop}}({\cal P})(a)) = {\sf tr}(r(a){\circ}{\cal P}{\circ}\ell (p)) = \epsilon (a)p(1) {\sf tr}{\cal E}({\cal P})
$$
for all $p \in A^*$ and $a \in A$ we conclude that ${\sf tr}{\cal E}({\cal P}) \neq 0$. In light of Lemma \ref{LemCalP} there are right integrals $\lambda \in \int^r$ and a $\Lambda \in \int_r$ such that $\lambda (\Lambda ) = 1$. This completes our proof of
\cite[Theorem 1]{RTrace}
which is based on the ideas of 
\cite[Section 5]{RTrace}
and does not use Hopf modules.

We conclude this section with a description of the endomorphism $P$ in terms of the variant ${\cal E}_{A^{cop}}$ of ${\cal E}$. Part c) of the following result is 
\cite[Lemma 3.3]{KuperRef}.
\begin{Theorem}\label{KIntCoInt}
Let $A$ be a finite-dimensional Hopf algebra over the field $k$ and suppose that $P$ is endomorphism of $A$ defined by (\ref{EqKupP}). Then:
\begin{enumerate}
\item[{\rm a)}]
$p(P(a)) = {\sf tr}(s{\circ}\ell (p){\circ}s{\circ}r(a)) = {\sf tr}(\ell (s(a)){\circ}s^2{\circ}\ell (p))$ for all $p \in A^*$ and $a \in A$.
\item[{\rm b)}]
$P = {\cal E}_{A^{cop}}(s^2){\circ}s$.
\item[{\rm c)}]
Integrals and cointegrals for $A$ are unique up to scalar multiple. In particular, if $\lambda \in \int^r$ and $\Lambda \in \int_r$ satisfy $\lambda (\Lambda ) = 1$ then $P = \lambda {\otimes}\Lambda$.
\end{enumerate}
\end{Theorem}

\pf
Since $s{\circ}r(a) = \ell (s(a)){\circ}s$ for all $a \in A$, to show part a) we need only establish the first equation. It is an easy exercise to show that
$$
(s{\circ}\ell (\alpha_\imath ){\circ}s{\circ}r(a_\jmath))(a_w) = m^{\ell}_{w, \jmath}s^u_{\ell}\Delta^{v, \imath}_us^x_va_x
$$
for all $1 \leq \imath, \jmath, w \leq n$. Therefore 
$$
\alpha_\imath (P(a_\jmath )) = P^\jmath_\imath = m^{\ell}_{w, \jmath}s^u_{\ell}\Delta^{v, \imath}_us^w_v = {\sf tr}(s{\circ}\ell (\alpha_\imath ){\circ}s{\circ}r(a_\jmath ))
$$ 
for all $1 \leq \imath, \jmath \leq n$ by (\ref{EqKupP}) which establishes the first equation of part a). Part b) follows from part a).

Since ${\sf tr}P = 1$, to show part c) we need only show that any integral and cointegral $P' = {\cal E}(f)$ for $A$ is a scalar multiple of ${\cal E}(s^{-2})$, or equivalently that $f = \alpha s^{-2}$ for some $\alpha \in k$.  By parts b)--d) of Proposition \ref{EProperties} the conditions $f{\bullet}a = \epsilon (a)f$ for all $a \in A$ and $f{\bullet}p = p(1)f$ for all $p \in A^*$ hold for $f$. Since $f(s^2(a)) = ((f{\bullet}s(a_{(1)}))(1))a_{(2)}$ for all $a \in A$ holds for any endomorphism of $A$, as a consequence of the first condition $f = \ell (f(1)){\circ}s^{-2}$. Since $(f{\bullet}p)^* = f^*{\bullet} p$ for all $p \in A^*$, the second condition therefore implies $f^* = r(\epsilon {\circ}f)^*{\circ}S^{-2} = (s^{-2}{\circ}r(\epsilon {\circ}f))^*$, or equivalently that $f = s^{-2}{\circ}r(\epsilon {\circ}f)$. By the last equation $f(1) = \epsilon (f(1))1$ which means $f = \ell (f(1)){\circ}s^{-2} = \alpha s^{-2}$, where $\alpha = \epsilon (f(1))$.
\qed

\end{document}